\documentclass[12pt]{article}

\usepackage{amsmath, epsfig, cite}
\usepackage{amssymb}
\usepackage{amsfonts}
\usepackage{latexsym}
\usepackage{float}
\usepackage{color}

\newtheorem{thm}{Theorem}[section]

\newtheorem{lem}[thm]{Lemma}

\numberwithin{equation}{section}

\newcommand{\qed}{{\hfill$\square$}\medskip}
\begin{document}

\begin{center}
{\large\bf Further results on the divisibility of $q$-trinomial coefficients}
\end{center}

\vskip 2mm \centerline{Ji-Cai Liu and Wei-Wei Qi}
\begin{center}
{\footnotesize Department of Mathematics, Wenzhou University, Wenzhou 325035, PR China\\[5pt]
{\tt jcliu2016@gmail.com, 2484466171@qq.com} \\[10pt]
}
\end{center}

\vskip 0.7cm \noindent{\bf Abstract.}
We study divisibility for the $q$-trinomial coefficients $\tau_0(n,m,q)$, $T_0(n,m,q)$
and $T_1(n,m,q)$, which were first introduced by Andrews and Baxter. In particular, we completely determine $\tau_0(an,bn,q)$, $T_0(an,bn,q)$ and $T_1(an,bn,q)$ modulo the square of the cyclotomic polynomial $\Phi_n(q)$ for $(a,b)=(m,m-1)$.


\vskip 2mm
\noindent{\it MR Subject Classifications}: 11A07, 11B65, 13A05, 05A10

\section{Introduction}
In 1987, Andrews and Baxter \cite{ab-1987} introduced six kinds of $q$-trinomial coefficients in the study of the solution of a model in statistical mechanics, which can be listed as follows:
\begin{align*}
&\left(\!\!\binom n m \!\!\right)_q=\sum_{k=0}^{n}
q^{k(k+m)}{n \brack k} { n-k \brack k+m},\\[10pt]
&\tau_0(n,m,q)=\sum_{k=0}^n(-1)^kq^{nk-{k\choose 2}}{n\brack k}{2n-2k\brack n-m-k},\\[10pt]
&T_0(n,m,q)=\sum_{k=0}^n(-1)^k{n\brack k}_{q^2}{2n-2k\brack n-m-k},\\[10pt]
&T_1(n,m,q)=\sum_{k=0}^n(-q)^k{n\brack k}_{q^2}{2n-2k\brack n-m-k},\\[10pt]
&t_0(n,m,q)=\sum_{k=0}^n(-1)^kq^{k^2}{n\brack k}_{q^2}{2n-2k\brack n-m-k},\\[10pt]
&t_1(n,m,q)=\sum_{k=0}^n(-1)^kq^{k(k-1)}{n\brack k}_{q^2}{2n-2k\brack n-m-k}.
\end{align*}
Here and in what follows, the $q$-binomial coefficients are defined as
$$
{n\brack k}={n\brack k}_q
=\begin{cases}
\displaystyle\frac{(1-q^n)(1-q^{n-1})\cdots (1-q^{n-k+1})}{(1-q)(1-q^2)\cdots(1-q^k)}, &\text{if $0\leqslant k\leqslant n$},\\[10pt]
0,&\text{otherwise.}
\end{cases}
$$

Note that these $q$-trinomial coefficients are six apparently distinct $q$-analogues of the trinomial coefficients $\left(\!\binom n m \!\right)$, which are given by
\begin{align*}
(1+x+x^2)^n=\sum_{m=-n}^{n}\left(\!\!\binom n m \!\!\right)x^{m+n}.
\end{align*}
It is well-known that the trinomial coefficients possess the following two simple formulas (see \cite[page 43]{sills-b-2018}):
\begin{align*}
\left(\!\!\binom n m \!\!\right)=\sum_{k=0}^n{n\choose k}{n-k\choose k+m},
\end{align*}
and
\begin{align*}
\left(\!\!\binom n m \!\!\right)=\sum_{k=0}^n(-1)^k{n\choose k}{2n-2k\choose n-m-k}.
\end{align*}

In the past two decades, $q$-analogues of congruences ($q$-congruences) were widely studied by many researchers. For recent developments on $q$-congruences, we refer the interested reader to
\cite{bachraoui-rj-2021,goro-ijnt-2019,guo-aam-2020,gs-bams-2022,gz-am-2019,
lw-2020,lj-2021,lw-jmaa-2021,lw-rm-2022,ni-rmjm-2021,ps-aam-2021,wang-rm-2021,WY2,wei-jcta-2021}.

It is remarkable that Andrews \cite{andrews-dm-1999} showed that for any odd prime $p$,
\begin{align}
{2p-1\brack p-1}\equiv q^{\frac{p(p-1)}{2}}\pmod{[p]_q^2},\label{newa-9}
\end{align}
which gives a $q$-analogue of Babbage's congruence \cite{babbage-epj-1819}.
In order to understand \eqref{newa-9}, we recall some necessary notation.
For polynomials $A_1(q), A_2(q),P(q)\in \mathbb{Z}[q]$, the $q$-congruence $$A_1(q)/A_2(q)\equiv 0\pmod{P(q)}$$ is understood as $A_1(q)$ is divisible by $P(q)$ and $A_2(q)$ is coprime with $P(q)$. In general, for rational functions $A(q),B(q)\in \mathbb{Z}(q)$,
\begin{align*}
A(q)\equiv B(q)\pmod{P(q)}\Longleftrightarrow
A(q)-B(q)\equiv 0\pmod{P(q)}.
\end{align*}
The $q$-integers are defined as $[n]_q=(1-q^n)/(1-q)$ for $n\ge 1$, and the $n$th
cyclotomic polynomial is given by
\begin{align*}
\Phi_n(q)=\prod_{\substack{1\le k \le n\\[3pt](n,k)=1}}
(q-\zeta^k),
\end{align*}
where $\zeta$ denotes an $n$th primitive root of unity.

It is worth mentioning that Straub \cite[Theorem 2.2]{straub-pams-2019} extended \eqref{newa-9} as follows (notice that ${2n-1\brack n-1}={2n\brack n}/(1+q^n)$):
\begin{align}
{an\brack bn}\equiv {a\brack b}_{q^{n^2}}-(a-b)b{a\choose b}\frac{n^2-1}{24}(q^n-1)^2\pmod{\Phi_n(q)^3},
\label{newa-1}
\end{align}
which was further generalized by Zudilin \cite{zudilin-ac-2019}.

The first author \cite{liu-rj-2022} investigated congruence properties for the
$q$-trinomial coefficients $\left(\!\binom {an} {bn} \!\right)_q$ for $(a,b)\in \{(1,0),(2,1)\}$ and showed that for any positive integer $n$,
\begin{align}
\left(\!\!\binom n 0 \!\!\right)_q\equiv \mathcal{A}_n(q)\pmod{\Phi_n(q)^2},\label{aaa-1}
\end{align}
and
\begin{align}
\left(\!\!\binom {2n} n \!\!\right)_q\equiv
2\mathcal{A}_n(q)-n(1-q^n)\pmod{\Phi_n(q)^2},\label{aaa-2}
\end{align}
where $\mathcal{A}_n(q)$ is given by
\begin{align*}
\mathcal{A}_n(q)={\begin{cases}
\displaystyle (-1)^m(1+q^m)q^{m(3m-1)/2}, &\text{if}~~n=3m,\\[10pt]
(-1)^mq^{m(3m+1)/2}, &\text{if}~~n=3m+1,\\[10pt]
(-1)^{m+1}q^{(m+1)(3m+2)/2},&\text{if}~~n=3m+2.
\end{cases}}
\end{align*}
It is remarkable that Chen, Xu and Wang \cite{cxw-2022} completely determined $\left(\!\binom {mn} {(m-1)n} \!\right)_q$ modulo $\Phi_n(q)^2$, which includes \eqref{aaa-1} and \eqref{aaa-2} as special cases.

Recently, the first author, Zhao and Jiang \cite{lzj-2022} further studied congruences for
$\tau_0(an,bn,q)$, $T_0(an,bn,q)$ and $T_1(an,bn,q)$ modulo $\Phi_n(q)^2$ for $(a,b)\in \{(1,0),(2,1)\}$.

The motivation of the paper is to completely determine $\tau_0(an,bn,q)$, $T_0(an,bn,q)$ and $T_1(an,bn,q)$ modulo $\Phi_n(q)^2$ for $(a,b)=(m,m-1)$, which extends the first author, Zhao and Jiang's results. The main results consist of the following three theorems.

\begin{thm}\label{t-1}
If $m$ and $n$ are both positive integers, then the following holds $\Phi_n(q)^2$:
\begin{align*}
\tau_0(mn,(m-1)n,q)&\equiv 2m-nm(2m-1)(1-q^n)\\[10pt]
&+(-1)^nq^{n(n+1)/2}\left(m(\mathcal{A}_n(q)+\mathcal{B}_n(q)-1)-\frac{3nm(m-1)}{2}(1-q^n)\right),
\end{align*}
where $\mathcal{B}_n(q)$ is given by
\begin{align*}
\mathcal{B}_n(q)={\begin{cases}
\displaystyle (-1)^m(1+q^{2m})q^{m(3m-5)/2}, &\text{if}~~n=3m,\\[10pt]
(-1)^mq^{m(3m+1)/2}, &\text{if}~~n=3m+1,\\[10pt]
(-1)^{m+1}q^{(m-1)(3m+2)/2},&\text{if}~~n=3m+2.
\end{cases}}
\end{align*}
\end{thm}

\begin{thm}\label{t-2}
If $m$ and $n$ are both positive integers, then the following holds $\Phi_n(q)^2$:
\begin{align*}
&T_0(mn,(m-1)n,q)\equiv
2m-nm(2m-1)(1-q^n)\\[10pt]
&+(-1)^n\left((1+(-1)^n)(m-1)+1\right)\left(m-nm(m-1)(1-q^n)+2m(\mathcal{A}_n(q)-1)\right).
\end{align*}
\end{thm}

\begin{thm}\label{t-3}
If $m$ and $n$ are both positive integers, then the following holds $\Phi_n(q)^2$:
\begin{align*}
&T_1(mn,(m-1)n,q)\equiv
2m-nm(2m-1)(1-q^n)\\[10pt]
&+(-1)^n\left((1+(-1)^n)(m-1)+1\right)\left(mq^n-nm(m-1)(1-q^n)+2m(\mathcal{B}_n(q)-1)\right).
\end{align*}
\end{thm}

The rest of the paper is organized as follows. In Section 2, we first establish some preliminary results.
The proofs of Theorems \ref{t-1}--\ref{t-3} will be given in Sections 3--5, respectively.

\section{Preliminary results}
In order to prove Theorems \ref{t-1}--\ref{t-3}, we first require two $q$-binomial identities.

\begin{lem} (See \cite[Lemma 2.3]{liu-cmj-2017}.)
For any non-negative integer $n$, we have
\begin{align}
&(1-q^n)\sum_{k=0}^{\lfloor n/2\rfloor}\frac{(-1)^kq^{k(k-1)/2}}{1-q^{n-k}}{n-k\brack k}\notag\\[10pt]
&={\begin{cases}
\displaystyle (-1)^m(1+q^m)q^{m(3m-1)/2}, &\text{if}~~n=3m,\\[10pt]
(-1)^mq^{m(3m+1)/2}, &\text{if}~~n=3m+1,\\[10pt]
(-1)^{m+1}q^{(m+1)(3m+2)/2},&\text{if}~~n=3m+2.
\end{cases}}\label{newb-6}
\end{align}
\end{lem}

\begin{lem} (See \cite[Lemma 2.3]{lzj-2022}.)
For any non-negative integer $n$, we have
\begin{align}
&(1-q^n)\sum_{k=0}^{\lfloor n/2\rfloor}\frac{(-1)^kq^{k(k-3)/2}}{1-q^{n-k}}{n-k\brack k}\notag\\[10pt]
&={\begin{cases}
\displaystyle (-1)^m(1+q^{2m})q^{m(3m-5)/2}, &\text{if}~~n=3m,\\[10pt]
(-1)^mq^{m(3m+1)/2}, &\text{if}~~n=3m+1,\\[10pt]
(-1)^{m+1}q^{(m-1)(3m+2)/2},&\text{if}~~n=3m+2.
\end{cases}}\label{newb-3}
\end{align}
\end{lem}

We remark that an induction proof of \eqref{newb-3} was presented in \cite{lzj-2022}, and
Chu \cite[Theorem 5]{chu-bmmss-2022} recently gave a common generalization of \eqref{newb-6} and \eqref{newb-3} through generating function technique.
Note that $\mathcal{A}_n(q)$ and $\mathcal{B}_n(q)$ coincide with the right-hand sides of \eqref{newb-6} and \eqref{newb-3}, respectively.

We also need the following congruence regarding $q$-binomial coefficients.
\begin{lem}
For positive integers $m$ and $n$, we have
\begin{align}
{mn\brack n}_{q^2}\equiv \left((1+(-1)^n)(m-1)+1\right)\left(m-nm(m-1)(1-q^n)\right)\pmod{\Phi_n(q)^2}.
\label{bbb-2}
\end{align}
\end{lem}
{\noindent\it Proof.}
By \cite[(2.7)]{cxw-2022}, we have
\begin{align}
{mn\brack n}\equiv m-\frac{nm(m-1)}{2}(1-q^n)\pmod{\Phi_n(q)^2}.\label{bbb-1}
\end{align}
It is clear that
\begin{align}
q^n\equiv 1\pmod{\Phi_n(q)}.\label{bbb-4}
\end{align}

If $n$ is odd, then $\Phi_n(q)|\Phi_n(q^2)$.
It follows from \eqref{bbb-1} and \eqref{bbb-4} that
\begin{align*}
{mn\brack n}_{q^2}&\equiv m-\frac{nm(m-1)}{2}(1+q^n)(1-q^n)\\[10pt]
&\equiv m-nm(m-1)(1-q^n)  \pmod{\Phi_n(q)^2},
\end{align*}
which proves the case $n\equiv 1\pmod{2}$ of \eqref{bbb-2}.

From \eqref{newa-1} and the fact $\Phi_{2n}(q)|\Phi_n(q^2)$, we deduce that
\begin{align}
{2mn\brack 2n}_{q^2}\equiv {2m\brack 2}_{q^{2n^2}}\pmod{\Phi_{2n}(q)^2}.\label{newb-8}
\end{align}
Note that for any positive integer $s$,
\begin{align*}
q^{2sn^2}&=1-\left(1-(q^{2n})^{sn}\right)\notag\\[10pt]
&=1-(1-q^{2n})(1+q^{2n}+q^{4n}+\cdots+q^{2n(sn-1)})\notag\\[10pt]
&\equiv 1-sn(1-q^{2n})\pmod{\Phi_{2n}(q)^2}.
\end{align*}
Thus,
\begin{align}
{2m\brack 2}_{q^{2n^2}}&=\sum_{i=0}^{2m-2}q^{2in^2}\sum_{j=0}^{m-1}q^{4jn^2}\notag\\[10pt]
&\equiv \left(2m-1-n(2m-1)(m-1)(1-q^{2n})\right)\left(m-nm(m-1)(1-q^{2n})\right)\notag\\[10pt]
&\equiv (2m-1)\left(m-2nm(m-1)(1-q^{2n})\right)\pmod{\Phi_{2n}(q)^2}.\label{bbb-3}
\end{align}
It follows from \eqref{newb-8} and \eqref{bbb-3} that for even positive integer $n$,
\begin{align*}
{mn\brack n}_{q^2}\equiv (2m-1)\left(m-nm(m-1)(1-q^{n})\right)\pmod{\Phi_{n}(q)^2},
\end{align*}
which is the case $n\equiv 0\pmod{2}$ of \eqref{bbb-2}.
\qed

\section{Proof of Theorem \ref{t-1}}
Note that
\begin{align}
&\tau_0(mn,(m-1)n,q)\notag\\[10pt]
&=\sum_{k=0}^{n}(-1)^kq^{mnk-{k\choose 2}}{mn\brack k}{2mn-2k\brack n-k}\notag\\[10pt]
&=(-1)^{n}\sum_{k=0}^{n}(-1)^kq^{(n-k)((2m-1)n+k+1)/2}{mn\brack n-k}{2n(m-1)+2k\brack k}\notag\\[10pt]
&=(-1)^nq^{n((2m-1)n+1)/2}{mn\brack n}+{2mn\brack n}\notag\\[10pt]
&+(-1)^n\sum_{k=1}^{n-1}(-1)^kq^{(n-k)((2m-1)n+k+1)/2}{mn\brack n-k}{2n(m-1)+2k\brack k},\label{newc-6}
\end{align}
where we have performed the variable substitution $k\to n-k$ in the second step.

For $1\le k\le n-1$, by \eqref{bbb-4} we have
\begin{align}
{mn\brack n-k}&={mn\brack n}\frac{(1-q^{n-k+1})(1-q^{n-k+2})\dots
(1-q^{n})}{(1-q^{(m-1)n+1})(1-q^{(m-1)n+2})\cdots(1-q^{(m-1)n+k})}\notag\\[10pt]
&\equiv{mn\brack n}\frac{(-1)^{k-1}q^{-k(k-1)/2}(1-q^n)}{1-q^k}\notag\\[10pt]
&\equiv{mn\brack n}\frac{(-1)^{k}q^{-k(k+1)/2}(1-q^n)}{1-q^{n-k}}\pmod{\Phi_n(q)^2}.\label{newc-8}
\end{align}
Furthermore, by \cite[Lemma 3.3]{tauraso-aam-2012} we have
\begin{align}
{2k-1\brack k}\equiv (-1)^kq^{k(3k-1)/2}{n-k\brack k} \pmod{\Phi_n(q)},\label{newc-4}
\end{align}
for $1\le k\le n-1$.
It follows from \eqref{bbb-4} and \eqref{newc-4} that for $1\le k\le n-1$,
\begin{align}
{2n(m-1)+2k\brack k}
&=\frac{(1-q^{2n(m-1)+k+1})(1-q^{2n(m-1)+k+2})\dots (1-q^{2n(m-1)+2k})}{(1-q)(1-q^{2})\cdots(1-q^{k})}\notag\\[10pt]
&\equiv(1+q^{k}){2k-1\brack k}\notag\\[10pt]
&\equiv (-1)^kq^{k(3k-1)/2}(1+q^{k}){n-k\brack k} \pmod{\Phi_n(q)}.\label{newc-9}
\end{align}
Combining \eqref{newc-8} and \eqref{newc-9} with the fact that
\begin{align*}
\frac{(n-k)((2m-1)n+k+1)}{2}=-\frac{k(k+1)}{2}-n(m-1)(k-n)+\frac{n(n+1)}{2},
\end{align*}
we arrive at
\begin{align}
&\sum_{k=1}^{n-1}(-1)^kq^{(n-k)((2m-1)n+k+1)/2}{mn\brack n-k}{2n(m-1)+2k\brack k}\notag\\[10pt]
&\equiv q^{n(n+1)/2}(1-q^n){mn\brack n}\sum_{k=1}^{n-1}\frac{(-1)^kq^{k(k-3)/2}(1+q^{k})}{1-q^{n-k}}{n-k\brack k} \notag\\[10pt]
&=q^{n(n+1)/2}{mn\brack n}\left(\mathcal{A}_n(q)+\mathcal{B}_n(q)-2\right)\pmod{\Phi_n(q)^2},\label{newc-10}
\end{align}
where we have used \eqref{newb-6} and \eqref{newb-3} in the last step.

Noting that
\begin{align*}
q^{n((2m-1)n+1)/2}=q^{n(n+1)/2+(m-1)n^2},
\end{align*}
and
\begin{align*}
q^{(m-1)n^2}&=1-(1-q^{(m-1)n})(1+q^{(m-1)n}+q^{2(m-1)n}+\cdots+q^{(m-1)n(n-1)})\\[10pt]
&\equiv 1-n(1-q^{(m-1)n})\\[10pt]
&= 1-n(1-q^n)(1+q^n+q^{2n}+\cdots+q^{(m-2)n})\\[10pt]
&\equiv 1-n(m-1)(1-q^n)\pmod{\Phi_n(q)^2},
\end{align*}
we obtain
\begin{align}
q^{n((2m-1)n+1)/2}\equiv q^{n(n+1)/2}\left(1-n(m-1)(1-q^n)\right)\pmod{\Phi_n(q)^2}.
\label{bbb-5}
\end{align}

From \eqref{newb-6} and \eqref{newb-3}, we deduce that
\begin{align}
\mathcal{A}_n(q)-1\equiv 0\pmod{\Phi_n(q)},\label{bbb-9}
\end{align}
and
\begin{align}
\mathcal{B}_n(q)-1\equiv 0\pmod{\Phi_n(q)},\label{bbb-10}
\end{align}
and so
\begin{align}
\mathcal{A}_n(q)+\mathcal{B}_n(q)-2\equiv 0\pmod{\Phi_n(q)}.\label{bbb-6}
\end{align}

Finally, substituting \eqref{bbb-1}, \eqref{newc-10} and \eqref{bbb-5} into the right-hand side of
\eqref{newc-6} and using \eqref{bbb-6}, we complete the proof of Theorem \ref{t-1}.

\section{Proof of Theorem \ref{t-2}}
Note that
\begin{align}
&T_0(mn,(m-1)n,q)\notag\\[10pt]
&=\sum_{k=0}^{n}(-1)^k{mn\brack k}_{q^2}{2mn-2k\brack n-k}\notag\\[10pt]
&=(-1)^n\sum_{k=0}^{n}(-1)^k{mn\brack n-k}_{q^2}{2n(m-1)+2k\brack k}\notag\\[10pt]
&=(-1)^n{mn\brack n}_{q^2}+{2mn\brack n}+(-1)^n\sum_{k=1}^{n-1}(-1)^k{mn\brack n-k}_{q^2}{2n(m-1)+2k\brack k}.\label{newd-4}
\end{align}
For $1\le k \le n-1$, we have
\begin{align}
{mn\brack n-k}_{q^2}&={mn\brack n}_{q^2}\frac{(1-q^{2(n-k+1)})(1-q^{2(n-k+2)})\dots
(1-q^{2n})}{(1-q^{2(m-1)n+2})(1-q^{2(m-1)n+4})\cdots(1-q^{2(m-1)n+2k})}\notag\\[10pt]
&\equiv{mn\brack n}_{q^2}\frac{(-1)^{k-1}q^{-k(k-1)}(1-q^{2n})}{1-q^{2k}}\notag\\[10pt]
&\equiv {mn\brack n}_{q^2}\frac{2(-1)^{k}q^{-k^2}(1-q^n)}{(1+q^k)(1-q^{n-k})}\pmod{\Phi_n(q)^2}.
\label{bbb-7}
\end{align}
It follows from \eqref{newb-6}, \eqref{newc-9} and \eqref{bbb-7} that
\begin{align}
&\sum_{k=1}^{n-1}(-1)^k{mn\brack n-k}_{q^2}{2n(m-1)+2k\brack k}\notag\\[10pt]
&\equiv 2(1-q^n){mn\brack n}_{q^2}\sum_{k=1}^{n-1} \frac{(-1)^kq^{k(k-1)/2}}{1-q^{n-k}}
{n-k\brack k}\notag\\[10pt]
&=2{mn\brack n}_{q^2}\left(\mathcal{A}_n(q)-1\right)\pmod{\Phi_n(q)^2}.\label{bbb-8}
\end{align}

Finally, substituting \eqref{bbb-2}, \eqref{bbb-1} and \eqref{bbb-8} into the right-hand side of
\eqref{newd-4} and using \eqref{bbb-9}, we complete the proof of Theorem \ref{t-2}.

\section{Proof of Theorem \ref{t-3}}
Note that
\begin{align}
&T_1(mn,(m-1)n,q)\notag\\[10pt]
&=\sum_{k=0}^{n}(-q)^k{mn\brack k}_{q^2}{2mn-2k\brack n-k}\notag\\[10pt]
&=\sum_{k=0}^{n}(-q)^{n-k}{mn\brack n-k}_{q^2}{2n(m-1)+2k\brack k}\notag\\[10pt]
&=(-q)^n{mn\brack n}_{q^2}+{2mn\brack n}+\sum_{k=1}^{n-1}(-q)^{n-k}{mn\brack n-k}_{q^2}{2n(m-1)+2k\brack k}.\label{bbb-11}
\end{align}
Similarly to the proof of Theorem \ref{t-2}, by using \eqref{newb-3}, \eqref{bbb-7} and \eqref{newc-9} we
obtain
\begin{align}
&\sum_{k=1}^{n-1}(-q)^{n-k}{mn\brack n-k}_{q^2}{2n(m-1)+2k\brack k}\notag\\[10pt]
&\equiv  2(-1)^n{mn\brack n}_{q^2}\left(\mathcal{B}_n(q)-1\right) \pmod{\Phi_n(q)^2}.
\label{bbb-12}
\end{align}
Substituting \eqref{bbb-2}, \eqref{bbb-1} and \eqref{bbb-12} into the right-hand side of
\eqref{bbb-11} and using \eqref{bbb-10}, we complete the proof of Theorem \ref{t-3}.

\vskip 5mm \noindent{\bf Acknowledgments.}
The first author was supported by the National Natural Science Foundation of China (grant 12171370).

\end{document}